# ROBUST GLOBAL ADAPTIVE EXPONENTIAL STABILIZATION OF DISCRETE-TIME SYSTEMS WITH APPLICATION TO FREEWAY TRAFFIC CONTROL

**Iasson Karafyllis[*], Maria Kontorinaki[**] and Markos Papageorgiou[**]**

[*]Dept. of Mathematics, National Technical University of Athens,
Zografou Campus, 15780, Athens, Greece, email: iasonkar@central.ntua.gr

[**]Dynamic Systems and Simulation Laboratory, Technical University of Crete,
Chania, 73100, Greece, email: mkontorinaki@dssl.tuc.gr, markos@dssl.tuc.gr

**Abstract**

This paper is devoted to the development of adaptive control schemes for uncertain discrete-time systems, which guarantee robust, global, exponential convergence to the desired equilibrium point of the system. The proposed control scheme consists of a nominal feedback law, which achieves robust, global, exponential stability properties when the vector of the parameters is known, in conjunction with a nonlinear, dead-beat observer. The obtained results are applicable to highly nonlinear, uncertain discrete-time systems with unknown constant parameters. The applicability of the obtained results to real control problems is demonstrated by the rigorous application of the proposed adaptive control scheme to uncertain freeway models.

## 1. Introduction

Adaptive control for discrete-time systems has been studied in many works (see for instance [8,29,30,31]) and in many cases it is a direct extension of adaptive control schemes for continuous-time systems (see [16]). Although discrete-time systems allow a direct study of the limitations of adaptive control schemes (see for example [28]), the major shortcoming of many adaptive control methodologies is that the closed-loop system does not exhibit an exponential convergence rate to the desired equilibrium point of the system, even if the nominal feedback law achieves global exponential stability properties when the parameters are precisely known.

This work is devoted to the development of adaptive control schemes for uncertain discrete-time systems, which guarantee robust, global, exponential convergence to the desired equilibrium point of the system. The idea is simple: use a nominal feedback law, which achieves robust, global, exponential stability properties when the vector of the parameters is known, in conjunction with a nonlinear, dead-beat observer. The dead-beat observer (designed using an extension of the methodology described in [12]) achieves the precise knowledge of the vector of unknown parameters after a transient period; then the states of the closed-loop system are robustly led to the desired equilibrium point with an exponential rate by the nominal feedback law. The proposed adaptive scheme does not require the knowledge of a Lyapunov function for the closed-loop system under the action of the nominal feedback stabilizer.

The obtained results are applicable to highly nonlinear, uncertain discrete-time systems with unknown constant parameters. The applicability of the obtained results to real control problems is demonstrated by the rigorous application of the proposed adaptive control scheme to uncertain freeway models.

Traffic congestion in freeways leads to serious degradation of the infrastructure causing excessive delays, and impacting traffic safety and the environment. Extensive research has been conducted to investigate and develop traffic control measures which can tackle this phenomenon. However, measures such as ramp metering, variable speed limits or dynamic route guidance have to be driven by appropriate control strategies in order to achieve their target. Traffic control strategies such as nonlinear optimal control [2,5] and Model Predictive Control [1,9] have been extensively studied but they are highly demanding from the computational point of view. However, the efficiency of traffic operations can also be enhanced by explicit feedback control approaches such as the pioneering I-type regulator ALINEA [23] and its extensions [26,27], as well as other proposed feedback control algorithms in [10,11,24,25]. These explicit feedback control strategies should guarantee local stability properties for the desired uncongested equilibrium point of the freeway model.

A Lyapunov approach was adopted in [14,15], which led to the robust, global exponential stabilization of the uncongested equilibrium point of a nonlinear freeway model. The nonlinear freeway model considered in [14,15] is a generalization of various freeway models (see [4,5,19]), which are special cases of the model used in [14,15]. However, the nonlinear feedback stabilizer demands the knowledge of several model parameters, which are usually unknown. The present work proposes an adaptive control scheme, which is based on a dead-beat nonlinear observer and guarantees the robust, global exponential convergence rate to the desired uncongested point of the freeway model.

The structure of the present work is as follows. Section 2 is devoted to the development of the robust, global, exponential adaptive control scheme for nonlinear uncertain discrete-time systems. The obtained results are applied rigorously in Section 3 to uncertain freeway models for the robust, global, exponential attractivity of the (unknown) desired uncongested equilibrium point of the freeway model. The concluding remarks of the paper are given in Section 4.

**Notation.** Throughout this paper, we adopt the following notation:
* $\Re_+ := [0, +\infty)$. For every set $S$, $S^n = \underbrace{S \times ... \times S}_{n \text{ times}}$ for every positive integer $n$. $\Re_+^n := (\Re_+)^n$. For every $x \in \Re$, $[x]$ denotes the integer part of $x \in \Re$. For certain sets $S_1, S_2, ..., S_n$, the set $S_1 \times S_2 \times ... \times S_n$ is denoted by $\prod_{i=1}^n S_i$.
* Let $x, y \in \Re^n$. By $|x|$ we denote the Euclidean norm of $x \in \Re^n$ and by $x'$ we denote the transpose of $x \in \Re^n$.
* When $R$ is an index set, then by $(x_i; i \in R)$ we denote a vector with components all $x_i \in \Re$ with $i \in R$, in increasing order. For example, if $R = \{2, 5, 10\}$, then $(x_i; i \in R) = (x_2, x_5, x_{10})'$.

## 2. Exponential Stabilization Of Systems With Unknown Parameters

Consider the discrete-time system:
$$z^+ = F(d,z), z \in X \subseteq \Re^n, d \in D \qquad (2.1)$$
where $X \subseteq \Re^n$ is a non-empty closed set with $z^* \in X$, $D \subseteq \Re^l$ is a non-empty set, $F: D \times X \to X$ is a locally bounded mapping with $F(d, z^*) = z^*$ for all $d \in D$. In this work we adopt the following robust exponential stability notion (see similar notions in [7,13,18]).

**Definition 2.2:** *We say that $z^* \in X$ is Robustly Globally Exponentially Stable (RGES) for system (2.1) if there exist constants $M, \sigma > 0$ such that for every $z_0 \in X$, $\{d_i \in D\}_{i=0}^\infty$, the solution $z(t)$ of (2.1) with $z(0) = z_0$ corresponding to $\{d_i \in D\}_{i=0}^\infty$ satisfies $|z(t) - z^*| \le M \exp(-\sigma t) |z_0 - z^*|$ for all $t \ge 0$.*

We next consider discrete-time systems with uncertain constant parameters and outputs. Consider the discrete-time system:
$$x^+ = f(d, \theta, x, u), \ x \in S, d \in D, u \in U \qquad (2.2)$$
where $S \subseteq \Re^n$, $D \subseteq \Re^l$, $U \subseteq \Re^m$, $\Theta \subseteq \Re^q$ are non-empty sets and $f: D \times \Theta \times S \times U \to S$ is a locally bounded mapping. In this setting, $x \in S$ denotes the state of system (2.2), $d \in D$ is an unknown, time-varying input, $u \in U$ is the control input and $\theta \in \Theta$ denotes the vector of unknown, constant parameters. The measured output of the system is given by
$$y(t) = h(d(t), \theta, x(t)) \qquad (2.3)$$
where $h: D \times \Theta \times S \to \Re^k$ is a locally bounded mapping. We assume that $x^* \in S$ is an equilibrium point for system (2.2) and $d \in D$ is a vanishing perturbation, i.e., there exist vectors $y^* \in h(D \times \{\theta\} \times S)$ such that $f(d, \theta, x^*, u^*) = x^*$, $y^* = h(d, \theta, x^*)$ for all $d \in D$. Moreover, let $Y \subseteq \Re^k$ be a set with $h(D \times \Theta \times S) \subseteq Y$.

In what follows we denote by $y^{(p)}(t) = (y(t-1), y(t-2), \ldots, y(t-p))$ for certain positive integer $p > 0$ the "$p$-history" of the signal $y(t)$ (defined for all $t \ge p$).
The main result of this section provides sufficient conditions for dynamic, robust, global, exponential stabilization of the equilibrium point $x^* \in S$. The stabilizer is constructed under the following assumptions.

**(H1)** *Suppose that there exists a mapping $k: \Theta \times Y \to U$ such that $x^* \in S$ is RGES for the closed-loop system (2.2), (2.3) with $u = k(\theta, y)$.*

**(H2)** *Suppose that there exist a positive integer $p > 0$, a mapping $\Psi: Y \times A \to \Theta$ and a set $A \subseteq Y^p$ which contains all $w \in Y^p$ in a neighborhood of $(y^*, \ldots, y^*)$, such that for every sequence $\{(d(t), \hat{\theta}(t)) \in D \times \Theta\}_{t=0}^\infty$ and for every $x_0 \in S$, the solution $x(t)$ of (2.2), (2.3) with $u = k(\hat{\theta}, y)$, initial condition $x(0) = x_0$ corresponding to inputs $\{(d(t), \hat{\theta}(t)) \in D \times \Theta\}_{t=0}^\infty$ satisfies $\theta = \Psi(y(t), y^{(p)}(t))$ for all $t \ge p$ with $y^{(p)}(t) \in A$.*

**(H3)** *There exists a positive integer $m > 0$, such that for every sequence $\{(d(t), \hat{\theta}(t)) \in D \times \Theta\}_{t=0}^\infty$ and for every $x_0 \in S$, the solution $x(t)$ of (2.2), (2.3) with $u = k(\hat{\theta}, y)$, initial condition $x(0) = x_0$ corresponding to inputs $\{(d(t), \hat{\theta}(t)) \in D \times \Theta\}_{t=0}^\infty$ satisfies $y^{(p)}(t - i(t)) \in A$ for some $i(t) \in \{0, 1, \ldots, m\}$ and for all $t \ge m + p$.*

Assumption (H1) is a standard assumption, which guarantees the existence of a robust global exponential stabilizer when the vector of the parameters $\theta \in \Theta$ is known. Assumptions (H2)-(H3) are equivalent to complete, robust observability assumption of $\theta$ from the output given by (2.3) (see also [12]).

**Theorem 2.1:** *Consider system (2.2) with output given by (2.3) under assumption (H1), (H2), (H3). Moreover, suppose that the sets $f(D \times \Theta \times S \times U), Y, \Theta$ are bounded. Finally, assume that there exist a constant $L \geq 0$, neighborhoods $N_1 \subseteq \Re^n$ of $x^*$, $N_2 \subseteq \Re^k$ of $y^*$, $N_3 \subseteq \Re^q$ of $\theta$, such that the inequalities*

$$\left| f(d, \theta, x, k(\hat{\theta}, h(d, \theta, x))) - x^* \right| + \left| h(d, \theta, x) - y^* \right| \leq L \left| x - x^* \right| + L \left| \hat{\theta} - \theta \right|,$$

$$\left| \Psi(h(d, \theta, x), w) - \theta \right| \leq L \left| x - x^* \right| + L \sum_{i=1}^{p} \left| w_i - y^* \right|$$

*hold for all $x \in N_1 \cap S$, $d \in D$, $\hat{\theta} \in N_3 \cap \Theta$, $w_i \in N_2 \cap Y$ ($i = 1, \ldots, p$) with $w = (w_1, \ldots, w_p)$. Then the dynamic feedback stabilizer*

$$\begin{aligned} w_1^+ &= y \\ w_2^+ &= w_1 \\ &\vdots \\ w_p^+ &= w_{p-1} \\ \hat{\theta}^+ &= \begin{cases} \hat{\theta} & \text{if } w \notin A \\ \Psi(y, w) & \text{if } w \in A \end{cases} \\ u &= k(\hat{\theta}, y) \end{aligned} \quad (2.4)$$

*where $w = (w_1, \ldots, w_p) \in Y^p$, $\hat{\theta} \in \Theta$ achieves the following:*

1) *There exist constants $M, \sigma > 0$ such that for every sequence $\{d(i) \in D\}_{i=0}^{\infty}$ and for every $(x_0, w_0, \hat{\theta}_0) \in S \times Y^p \times \Theta$, the solution $(x(t), w(t), \hat{\theta}(t))$ of the closed-loop system (2.2), (2.3) with (2.4), initial condition $(x(0), w(0), \hat{\theta}(0)) = (x_0, w_0, \hat{\theta}_0)$ corresponding to input $\{d(i) \in D\}_{i=0}^{\infty}$ satisfies*

$$\begin{aligned} & \left| x(t) - x^* \right| + \sum_{i=1}^{p} \left| w_i(t) - y^* \right| + \left| \hat{\theta}(t) - \theta \right| \\ & \leq M \exp(-\sigma t) \left( \left| x(0) - x^* \right| + \sum_{i=1}^{p} \left| w_i(0) - y^* \right| + \left| \hat{\theta}(0) - \theta \right| \right) \end{aligned}, \text{ for all } t \geq 0 \quad (2.5)$$

2) *For every sequence $\{d(i) \in D\}_{i=0}^{\infty}$ and for every $(x_0, w_0, \hat{\theta}_0) \in S \times Y^p \times \Theta$, the solution $(x(t), w(t), \hat{\theta}(t))$ of the closed-loop system (2.2), (2.3) with (2.4), initial condition $(x(0), w(0), \hat{\theta}(0)) = (x_0, w_0, \hat{\theta}_0)$ corresponding to input $\{d(i) \in D\}_{i=0}^{\infty}$ satisfies $\hat{\theta}(t) = \theta$, for all $t \geq m + p + 1$.*

**Remark:** The dynamic feedback stabilizer (2.4) achieves dead-beat estimation of the vector of unknown parameters $\theta \in \Theta$. More, specifically, the variable $\hat{\theta}$ provides an estimate of the vector of unknown parameters $\theta \in \Theta$. Due to the dead-beat estimation, the exponential convergence property for the closed-loop system is preserved, as estimate (2.5) shows.

The proof of Theorem 2.1 relies on the following technical lemma.

**Lemma 2.2:** *Consider system (2.1) and suppose that the following hold:*
**i)** *There exist constants $M, \sigma > 0$ such that for every $z_0 \in \Omega$, $\{d(i) \in D\}_{i=0}^{\infty}$ the solution $z(t)$ of (2.1) with initial condition $z(0) = z_0$ corresponding to input $\{d(i) \in D\}_{i=0}^{\infty}$ satisfies $\left| z(t) - z^* \right| \leq M \left| z_0 - z^* \right| \exp(-\sigma t)$ for all $t \geq 0$.*
**ii)** *There exists an integer $N \geq 1$ such that for every $z_0 \in X$, $\{d(i) \in D\}_{i=0}^{\infty}$ and $t \geq N$ there exists $i(t) \in \{0, 1, \ldots, N\}$ for which the solution $z(t)$ of (2.1) with initial condition $z(0) = z_0$ corresponding to input $\{d(i) \in D\}_{i=0}^{\infty}$ satisfies $z(t - i(t)) \in \Omega$.*

**iii)** *There exists a constant $L \geq 1$, such that the inequality $|F(d,z) - z^*| \leq L|z - z^*|$ holds for all $d \in D$ and for all $z \in X$ in a neighborhood of $z^*$.*

*Then $z^* \in X$ is RGES for the uncertain system (2.1).*

**Proof:** By virtue of assumption (iii), there exists $\delta > 0$ such that the inequality $|F(d,z) - z^*| \leq L|z - z^*|$ holds for all $d \in D$ and $z \in A := \{ y \in X : |y - z^*| < \delta \}$. Since $F : D \times X \to X$ is a bounded mapping, there exists a constant $R > 0$ which satisfies

$$\sup\{|F(d,z)| : z \in X, d \in D\} \leq R \tag{2.6}$$

It follows from (2.6) and the triangle inequality that the following inequality holds:

$$\sup\left\{\frac{|F(d,z) - z^*|}{|z - z^*|} : d \in D, z \in X \setminus A\right\} \leq$$
$$\delta^{-1} \sup\{|F(d,z) - z^*| : z \in X, d \in D\} \leq \delta^{-1}(R + |z^*|) \tag{2.7}$$

Combining (2.7) and the fact that $|F(d,z) - z^*| \leq L|z - z^*|$ holds for all $d \in D$ and for all $z \in A$, we get:

$$|F(d,z) - z^*| \leq \max\left(L, \delta^{-1}(R + |z^*|)\right)|z - z^*|, \text{ for all } (d,z) \in D \times X \tag{2.8}$$

Let $z_0 \in X$ be an arbitrary vector and let $\{d(i) \in D\}_{i=0}^{\infty}$ be an arbitrary sequence. Consider the solution $z(t)$ of $z^+ = F(d,z)$ with initial condition $z(0) = z_0$ corresponding to input $\{d(i) \in D\}_{i=0}^{\infty}$. By virtue of assumption (ii), there exists $i(N) \in \{0,1,...,N\}$ with $z(N - i(N)) \in \Omega$. By virtue of assumption (i), we get:

$$|z(t) - z^*| \leq M|z(k) - z^*|\exp(-\sigma(t-k)), \text{ for all } t \geq k, \text{ where } k = N - i(N). \tag{2.9}$$

Notice that $k \in \{0,1,...,N\}$. Using induction and (2.8), we get

$$|z(t) - z^*| \leq \widetilde{L}^t|z_0 - z^*|, \text{ for all } t \geq 0, \tag{2.10}$$

where $\widetilde{L} := \max\left(L, \delta^{-1}(R + |z^*|)\right) \geq 1$. Combining (2.9), (2.10) and the fact that $k \in \{0,1,...,N\}$, we obtain:

$$|z(t) - z^*| \leq M\widetilde{L}^N \exp(\sigma N)|z_0 - z^*|\exp(-\sigma t), \text{ for all } t \geq 0 \tag{2.11}$$

Noticing that assumption (iii) guarantees that $z^* = F(d, z^*)$, we conclude that estimate (2.11) implies that $z^* \in X$ is RGES for the uncertain system (2.1). The proof is complete. ◁

We are now ready to provide the proof of Theorem 2.1.

**Proof of Theorem 2.1:** Let $\Phi(x)$ be the (possibly empty) set of all $w = (w_1, ..., w_p) \in Y^p$ for which there exist $\xi \in S$, $(d(i), \hat{\theta}(i)) \in D \times \Theta$, $i = 0, ..., p-1$, such that the vectors $\bar{x}(i)$, $i = 0, ..., p$, defined by the recursive formula

$$\bar{x}(0) = \xi, \ \bar{x}(i+1) = f(d(i), \theta, \bar{x}(i), k(\hat{\theta}(i), h(d(i), \theta, \bar{x}(i)))), \text{ for } i = 0, ..., p-1, \tag{2.12}$$

satisfy $\bar{x}(p) = x$ and $w_{p-i} = h(d(i), \theta, \bar{x}(i))$ for $i = 0, ..., p-1$.

All assumptions of Lemma 2.2 hold with $X = S \times Y^p \times \Theta$, $z = (x, w, \hat{\theta})$, $\Omega = \{(x, w, \theta) : w \in \Phi(x), x \in S\}$, $N = m + p + 1$, $z^* = (x^*, y^*, ..., y^*, \theta)$ and

$$F(d,z) := \begin{bmatrix} f(d, \theta, x, k(\hat{\theta}, h(d, \theta, x))) \\ h(d, \theta, x) \\ w_1 \\ \vdots \\ w_{p-1} \\ g(h(d, \theta, x), w, \hat{\theta}) \end{bmatrix},$$

where $g(h(d,\theta,x),w,\hat{\theta}) := \begin{cases} \hat{\theta} & \text{if } w \notin A \\ \Psi(h(d,\theta,x),w) & \text{if } w \in A \end{cases}$. We show next that assumptions (i), (ii) of Lemma 2.2 are direct consequences of assumptions (H1), (H2), (H3).

Let $\{d(i) \in D\}_{i=0}^{\infty}$ be an arbitrary sequence and let $(x_0, w_0, \hat{\theta}_0) \in \Omega$ be an arbitrary vector with $\hat{\theta}_0 = \theta$. Consider the solution $(x(t), w(t), \hat{\theta}(t))$ of the closed-loop system (2.2), (2.3) with (2.4), initial condition $(x(0), w(0), \hat{\theta}(0)) = (x_0, w_0, \hat{\theta}_0)$ corresponding to input $\{d(i) \in D\}_{i=0}^{\infty}$. By virtue of (2.12), the component $x(t)$ of the solution satisfies $x(t) = \bar{x}(t+p)$ for all $t \geq 0$, for certain solution $\bar{x}(i)$ of the system $\bar{x}^+ = f(\delta, \theta, \bar{x}, k(v, h(d, \theta, \bar{x})))$ (that corresponds to certain inputs $\{(\delta(t), v(t)) \in D \times \Theta\}_{t=0}^{\infty}$ with $\delta(t+p) = d(t)$, $v(t+p) = \hat{\theta}(t)$ for all $t \geq 0$ and appropriate initial condition $\xi \in S$). Moreover, $w(t) = \bar{y}^{(p)}(t+p) \in \Phi(x(t))$ for all $t \geq 0$, where $\bar{y}(t) = h(\delta(t), \theta, \bar{x}(t))$. Notice that if $w(0) = w_0 \in A$ then $\bar{y}^{(p)}(p) \in A$, and, consequently, assumption (H2) guarantees that $\hat{\theta}(1) = \theta$. If $w(0) = w_0 \notin A$ then $\hat{\theta}(1) = \hat{\theta}(0) = \theta$. Using induction and the previous argument, it follows that $\hat{\theta}(t) = \theta$ for all $t \geq 0$. Therefore, assumption (i) of Lemma 2.1 is a consequence of assumption (H1).

Assumption (ii) of Lemma 2.1 follows from the fact that $w(t) = y^{(p)}(t) \in \Phi(x(t))$ for all $t \geq p$. Assumption (H3) guarantees that $w(t-i(t)) = y^{(p)}(t-i(t)) \in A$ for some $i(t) \in \{0,1,...,m\}$ and for all $t \geq m+p$. It follows from (2.4) that $\hat{\theta}(t-i(t)+1) = \theta$. Since $t-i(t)+1 \geq p+1$, we also get $w(t-i(t)+1) \in \Phi(x(t))$ and thus $z(t-i(t)+1) \in \Omega$. Therefore, assumption (ii) of Lemma 2.2 holds with $N = m+p+1$.

Since $A \subseteq Y^p$ contains all $w \in Y^p$ in a neighborhood of $(y^*,...,y^*)$ and since there exist neighborhoods $N_1 \subseteq \Re^n$ of $x^*$, $N_2 \subseteq \Re^k$ of $y^*$, $N_3 \subseteq \Re^q$ of $\theta$, such that the inequalities

$$\left|f(d,\theta,x,k(\hat{\theta},x))-x^*\right|+\left|h(d,\theta,x)-y^*\right| \leq L\left|x-x^*\right|+L\left|\hat{\theta}-\theta\right|, \quad \left|\Psi(h(d,\theta,x),w)-\theta\right| \leq L\left|x-x^*\right|+L\sum_{i=1}^{p}\left|w_i-y^*\right|$$

hold for all $x \in N_1 \cap S$, $d \in D$, $\hat{\theta} \in N_3 \cap \Theta$, $w_i \in N_2 \cap Y$ ($i=1,...,p$) with $w = (w_1,...,w_p)$, it follows that assumption (iii) of Lemma 2.2 holds. ◁

## 3. Application to Freeway Traffic Control

### *3.I. The freeway model*

We consider a freeway which consists of $n \geq 3$ components or cells; typical cell lengths may be 200-500 m. Each cell may have an external inflow (e.g. from corresponding on-ramps), located near the cell's upstream boundary; and an external outflow (e.g. via corresponding off-ramps), located near the cell's downstream boundary (Figure 1). The number of vehicles at time $t \geq 0$ in component $i \in \{1,...,n\}$ is denoted by $x_i(t)$. The total outflow and the total inflow of vehicles of the component $i \in \{1,...,n\}$ at time $t \geq 0$ are denoted by $F_{i,out}(t) \geq 0$ and $F_{i,in}(t) \geq 0$, respectively. All flows during a time interval are measured in [veh]. Consequently, the balance of vehicles (conservation equation) for each component $i \in \{1,...,n\}$ gives:

$$x_i(t+1) = x_i(t) - F_{i,out}(t) + F_{i,in}(t), \quad i=1,...,n, \quad t \geq 0. \tag{3.1}$$

Each component of the network has storage capacity $a_i > 0$ ($i=1,...,n$). Our first assumption states that the external (off-ramp) flows from each cell are constant percentages of the total exit flow, i.e., there exist constants $P_i \in [0,1)$, $i=1,...,n$, such that:

$$\begin{pmatrix} \text{flow of vehicles} \\ \text{from cell } i \text{ to cell } i+1 \end{pmatrix} = (1-P_i)F_{i,out}(t), \text{ for } i=1,...,n-1 \tag{3.2}$$

$$\begin{pmatrix} \text{flow of vehicles} \\ \text{from cell } i \text{ to regions out of the freeway} \end{pmatrix} = P_i F_{i,out}(t), \text{ for } i=1,...,n. \tag{3.3}$$

The constants $P_i$ are known as exit rates. Since the $n$-th cell is the last downstream cell of the considered freeway, we may assume that $P_n = 1$. We also assume that $P_i < 1$ for $i=1,...,n-1$, and that all exits to regions out of the network can accommodate the respective exit flows.

Our second assumption is dealing with the attempted outflows $f_i(x_i)$, i.e. the flows that will exit the cell if there is sufficient space in the downstream cell. We assume that there exist functions $f_i : [0, a_i] \to \Re_+$ with $0 < f_i(x_i) < x_i$ for $x_i \in (0, a_i]$, variables $s_i(t) \in [0,1]$, $i=2,...,n$, so that:

$$F_{i-1,out}(t) = s_i(t)f_{i-1}(x_{i-1}(t)), \; i=2,...,n, \; t \geq 0 \text{ and } F_{n,out}(t) = f_n(x_n(t)). \tag{3.4}$$

The variable $s_i(t) \in [0,1]$, for each $i=2,...,n$, indicates the percentage of the attempted outflow from cell $i-1$ that becomes actual outflow from the same cell. The function $f_i : [0, a_i] \to \Re_+$ is called, in the specialized literature of Traffic Engineering (see, e.g., [4,5,6,19,20,21]), the demand-part of the fundamental diagram of the $i$-th cell, i.e. the flow that will exit the cell $i$ if there is sufficient space in the downstream cell $i+1$. Notice that equation (3.4) for $F_{n,out}(t)$ follows from our assumption that all exits to regions out of the network can accommodate the exit flows.

Let $v_i \geq 0$ ($i=1,...,n$) denote the attempted external inflow to component $i \in \{1,...,n\}$ from the region out of the freeway. Typically, $v_i$, $i=2,...,n$, correspond to external on-ramp flows which may be determined by a ramp metering control strategy. For the very first cell 1, we assume, for convenience, that there is just one external inflow, $v_1 > 0$. Let the variables $W_i(t) \in [0,1]$, $i=1,...,n$, indicate the percentage of the attempted external inflow to component $i \in \{1,...,n\}$ that becomes actual inflow. Then, we obtain from (3.2) and (3.4):

$$F_{1,in}(t) = W_1(t)v_1(t) \text{ and } F_{i,in}(t) = W_i(t)v_i(t) + s_i(t)(1-P_{i-1})f_{i-1}(x_{i-1}(t)), \; i=2,...,n. \tag{3.5}$$

Our next assumption requires that the inflow of vehicles at the cell $i \in \{1,...,n\}$ at time $t \geq 0$, denoted by $F_{i,in}(t) \geq 0$, cannot exceed the supply function of cell $i \in \{1,...,n\}$ at time $t \geq 0$, i.e.,

$$F_{i,in}(t) \leq \min(q_i, c_i(a_i - x_i(t))), \; i=1,...,n, \; t \geq 0 \tag{3.6}$$

where $q_i \in (0,+\infty)$ denotes the maximum flow that the $i$-th cell can receive (or the capacity flow of the $i$-th cell) and $c_i \in (0,1]$ ($i=1,...,n$) denotes the congestion wave speed of the $i$-th cell.

Following [4], we assume that, when the total demand flow of a cell is lower than the supply of the downstream cell, i.e. when $v_i(t) + (1-P_{i-1})f_{i-1}(x_{i-1}(t)) \leq \min(q_i, c_i(a_i - x_i(t)))$ for some $i \in \{2,...,n\}$, then the demand flow can be fully accommodated by the downstream cell, and hence we have $s_i(t) = W_i(t) = 1$. Similarly, when $v_1(t) \leq \min(q_1, c_1(a_1 - x_1(t)))$, then we have $W_i(t) = 1$. In contrast, when the total demand flow of a cell is higher than the supply of the downstream cell, i.e. when $v_i(t) + (1-P_{i-1})f_{i-1}(x_{i-1}(t)) > \min(q_i, c_i(a_i - x_i(t)))$ for some $i \in \{2,...,n\}$ (or when $v_1(t) > \min(q_1, c_1(a_1 - x_1(t)))$), then the demand flow cannot be fully accommodated by the downstream cell, and the actual flow is determined by the supply function, i.e. we have $F_{i,in}(t) = \min(q_i, c_i(a_i - x_i(t)))$ (or $F_{1,in}(t) = \min(q_1, c_1(a_1 - x_1(t)))$). Therefore, we get:

$$F_{1,in}(t) = \min(q_1, c_1(a_1 - x_1(t)), v_1(t)), \; t \geq 0 \tag{3.7}$$

$$s_i(t) = (1 - d_i(t))\min\left(1, \max\left(0, \frac{\min(q_i, c_i(a_i - x_i(t))) - v_i(t)}{(1 - P_{i-1})f_{i-1}(x_{i-1}(t))}\right)\right) + d_i(t)\min\left(1, \frac{\min(q_i, c_i(a_i - x_i(t)))}{(1 - P_{i-1})f_{i-1}(x_{i-1}(t))}\right), \; i = 2,...,n, \; t \geq 0 \tag{3.8}$$

$$F_{i,in}(t) = \min(q_i, c_i(a_i - x_i(t)), v_i(t) + (1 - P_{i-1})f_{i-1}(x_{i-1}(t))), \; i = 2,...,n, \; t \geq 0 \tag{3.9}$$

where

$$d_i(t) \in [0,1], \; i = 2,...,n, \; t \geq 0 \tag{3.10}$$

are time-varying parameters. Note that, if the supply is higher than the total demand, then (3.8) yields $s_i = 1$, irrespective of the value of $d_i$, since the total demand flow can be accommodated by the downstream cell. Thus, the parameter $d_i$ determines the relative inflow priorities, when the downstream supply prevails. Specifically, when $d_i(t) = 0$, then the on-ramp inflow has absolute priority over the internal inflow; on the other hand, when $d_i(t) = 1$, then the internal inflow has absolute priority over the on-ramp inflow; while intermediate values of $d_i$ reflect intermediate priority cases. The parameters $d_i(t) \in [0,1]$ are treated as unknown parameters (disturbances). Notice that by introducing the parameters $d_i(t) \in [0,1]$ (and by allowing them to be time-varying), we have taken into account all possible cases for the relative priorities of the inflows (and we also allow the priority rules to be time-varying); see [3,4] for freeway models with specific priority rules, which are special cases of our general approach.

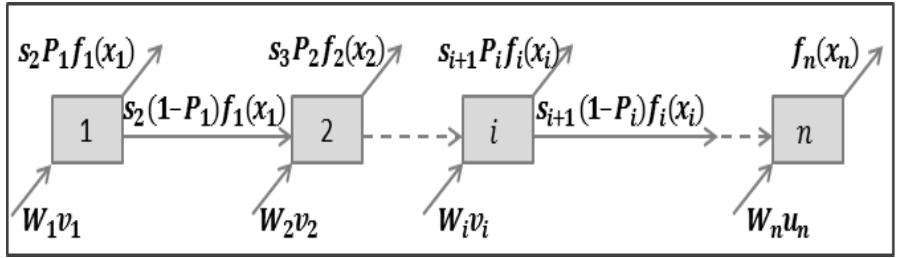

**Figure 1:** Scheme of the freeway model.

All the above are illustrated in Figure 1. Combining equations (3.1), (3.2), (3.3), (3.4), (3.7) and (3.9) we obtain the following discrete-time dynamical system:

$$\begin{aligned} x_1^+ &= x_1 - s_2 f_1(x_1) + \min(q_1, c_1(a_1 - x_1), v_1) \\ &= x_1 - s_2 f_1(x_1) + w_1 v_1 \end{aligned} \tag{3.11}$$

$$\begin{aligned} x_i^+ &= x_i - s_{i+1} f_i(x_i) + \min(q_i, c_i(a_i - x_i), v_i + (1 - P_{i-1})f_{i-1}(x_{i-1})) \\ &= x_i - s_{i+1} f_i(x_i) + W_i v_i + s_i(1 - P_{i-1})f_{i-1}(x_{i-1}) \end{aligned}, \text{ for } i = 2,...,n-1 \tag{3.12}$$

$$\begin{aligned} x_n^+ &= x_n - f_n(x_n) + \min(q_n, c_n(a_n - x_n), v_n + (1 - P_{n-1})f_{n-1}(x_{n-1})) \\ &= x_n - f_n(x_n) + W_n v_n + s_n(1 - P_{n-1})f_{n-1}(x_{n-1}) \end{aligned} \tag{3.13}$$

where $s_i \in [0,1]$, $i = 2,...,n$ are given by (3.8). The values of $W_i \in [0,1]$, $i = 1,...,n$, may also be similarly derived from (3.5), (3.7), (3.9) when $v_i > 0$ but they are not needed in what follows. Define $S = \prod_{i=1}^{n}(0, a_i]$. Since the functions $f_i : [0, a_i] \to \Re_+$ satisfy $0 < f_i(x_i) < x_i$ for all $x_i \in (0, a_i]$, it follows that (3.11), (3.12), (3.13) is an uncertain control system on $S$ (i.e., $x = (x_1,...,x_n)' \in S$) with inputs $v = (v_1,...,v_n)' \in (0,+\infty) \times \Re_+^{n-1}$ and disturbances $d = (d_2,...,d_n)' \in [0,1]^{n-1}$. We emphasize again that the uncertainty $d \in [0,1]^{n-1}$ appears in the equations (3.11), (3.12) and (3.13) only when the supply function prevails, i.e., only when $v_i(t) + (1 - P_{i-1})f_{i-1}(x_{i-1}(t)) > \min(q_i, c_i(a_i - x_i(t)))$ for some $i \in \{2,...,n\}$.

We make the following assumption for the functions $f_i : [0, a_i] \to \Re_+$ ($i = 1,...,n$):

**(H)** *There exist constants $\delta_i \in (0, a_i]$ and $r_i \in (0,1)$ such that $f_i(z) = r_i z$ for $z \in [0, \delta_i]$. Moreover, there exists a positive constant $f_i^{\min} > 0$ such that $f_i(\delta_i) = r_i \delta_i \geq f_i(z) \geq f_i^{\min}$ for all $z \in [\delta_i, a_i]$.*

Assumption (H) is a technical assumption that allows a very general class of demand functions (which are also allowed to be discontinuous). A more general assumption than assumption (H) was used in [14,15], but in [14,15] it was assumed that all parameters of the model were known. More specifically, in [14,15], it was not necessary the demand functions $f_i : [0, a_i] \to \Re_+$ ($i = 1,...,n$) to be linear on the corresponding intervals $[0, \delta_i]$.

*3.II. Global Exponential Stabilization of Freeway Models*

Define the vector field $\tilde{F} : D \times S \times (0, +\infty) \times \Re_+^{n-1} \to S$ for all $x \in S = \prod_{i=1}^{n}(0, a_i]$, $d = (d_2,...,d_n) \in D = [0,1]^{n-1}$ and $v \in (0, +\infty) \times \Re_+^{n-1}$:

$$\tilde{F}(d,x,v) = (\tilde{F}_1(d,x,v),...,\tilde{F}_n(d,x,v))' \in \Re^n$$

with $\tilde{F}_1(d,x,v) = x_1 - s_2 f_1(x_1) + \min(q_1, c_1(a_1 - x_1), v_1)$,

$\tilde{F}_i(d,x,v) = x_i - s_{i+1} f_i(x_i) + \min(q_i, c_i(a_i - x_i), v_i + (1 - P_{i-1}) f_{i-1}(x_{i-1}))$, for $i = 2,...,n-1$,

$\tilde{F}_n(d,x,v) = x_n - f_n(x_n) + \min(q_n, c_n(a_n - x_n), v_n + (1 - P_{n-1}) f_{n-1}(x_{n-1}))$ and

$$s_i = (1-d_i)\min\left(1, \max\left(0, \frac{\min(q_i, c_i(a_i - x_i)) - v_i}{(1 - P_{i-1}) f_{i-1}(x_{i-1})}\right)\right) + d_i \min\left(1, \frac{\min(q_i, c_i(a_i - x_i))}{(1 - P_{i-1}) f_{i-1}(x_{i-1})}\right), \text{ for } i = 2,...,n. \quad (3.14)$$

Notice that, using definition (3.14), the control system (3.11), (3.12), (3.13) can be written in the following vector form:

$$x^+ = \tilde{F}(d,x,v) \quad , \quad x \in S, d \in D, v \in (0, +\infty) \times \Re_+^{n-1}. \quad (3.15)$$

Consider the freeway model (3.15) under assumption (H). Let $v^* = (v_1^*,...,v_n^*)' \in (0, +\infty) \times \Re_+^{n-1}$ be a vector that satisfies:

$$v_1^* < \min(q_1, c_1(a_1 - \delta_1), r_1 \delta_1)$$

$$v_i^* + \sum_{j=1}^{i-1} v_j^* \left(\prod_{k=j}^{i-1}(1 - P_k)\right) < \min(q_i, c_i(a_i - \delta_i), r_i \delta_i) \quad (3.16)$$

Any inflow vector that satisfies (3.16), defines an uncongested equilibrium point $x^* = (x_1^*,...,x_n^*) \in \prod_{i=1}^{n}(0, \delta_i)$ for the freeway model:

$$x_1^* = r_1^{-1} v_1^*$$

$$x_i^* = r_i^{-1}\left(v_i^* + \sum_{j=1}^{i-1} v_j^* \left(\prod_{k=j}^{i-1}(1 - P_k)\right)\right) \quad , \quad i = 2,...,n \quad (3.17)$$

The uncongested equilibrium point is not globally exponentially stable for arbitrary $v_1^* > 0$, $v_i^* \geq 0$ ($i = 2,...,n$); indeed, for relatively large values of inflows $v_1^* > 0$, $v_i^* \geq 0$ ($i = 2,...,n$), other equilibria for model (3.15) (congested equilibria) may appear, for which the cell densities are large and can attract the solution of (3.15).

The following result (see [14,15]) is the main result in feedback design that provides the nominal feedback for the adaptive control scheme that we intend to use. The result shows that a continuous, robust, global exponential stabilizer exists for every freeway model of the form (3.15) under assumption (H).

**Theorem 3.1:** *Consider system (3.15) with $n \geq 3$ under assumption (H) for each $i = 1,...,n$. Then there exist a subset $R \subseteq \{1,...,n\}$ of the set of all indices $i \in \{1,...,n\}$ with $v_i^* > 0$, constants $\sigma \in (0,1]$, $b_i \in (0, v_i^*)$ for $i \in R$ and a constant $\tau^* > 0$ such that for every $\tau \in (0, \tau^*)$ the feedback law $k: S \to \Re_+^n$ defined by:*

$$k(x) = (k_1(x),...,k_n(x))' \in \Re^n \text{ with}$$

$$k_i(x) = \max\left(b_i, v_i^* - \tau^{-1}(v_i^* - b_i)\Xi(x)\right), \text{ for all } x \in S, i \in R \text{ and } k_i(x) = v_i^*, \text{ for all } x \in S, i \notin R \quad (3.18)$$

*where*

$$\Xi(x) := \sum_{i=1}^{n} \sigma^i \max\left(0, x_i - x_i^*\right), \text{ for all } x \in S \quad (3.19)$$

*achieves robust global exponential stabilization of the uncongested equilibrium point $x^*$ of system (3.15), i.e., $x^*$ is RGES for the closed-loop system (3.15) with $v = k(x)$.*

The result of Theorem 3.1 (see [14,15]) is based on the construction of a Control Lyapunov function for system (3.15) under a more general assumption than assumption (H). The feedback law provides values for the controllable inflows $v_i$, $i \in R$, in the interval $[b_i, v_i^*]$ for all $i \in R$, where $b_i \in (0, v_i^*)$ for $i \in R$ are the minimum allowable inflows. Since the proof of Theorem 3.1 is constructive, criteria for the selection of the index set $R \subseteq \{1,...,n\}$ and the constants $\sigma \in (0,1]$, $b_i \in (0, v_i^*)$ for $i \in R$ and $\tau^* > 0$ are provided.

Without loss of generality, we will assume, in what follows, that $R \neq \varnothing$ (because otherwise the uncongested equilibrium point is open-loop RGES).

Let $\mu_i \in (0, \delta_i)$, $v_{i,\max} \in (0, +\infty)$ ($i = 1,...,n$) be constants such that

$$v_{1,\max} < \min(q_1, c_1(a_1 - \mu_1)),$$
$$v_{i,\max} + (1 - P_{i-1})r_{i-1}\mu_{i-1} < \min(q_i, c_i(a_i - \mu_i)), \quad i = 2,...,n \quad (3.20)$$

It follows that if $x \in \Omega = \prod_{i=1}^{n}(0, \mu_i)$ and $v \in (0, v_{1,\max}] \times \prod_{i=2}^{n}[0, v_{i,\max}]$:

$$W_i = 1, \text{ for } i = 1,...,n \text{ and } s_i = 1, \text{ for } i = 2,...,n \quad (3.21)$$

$$x_1^+ = x_1 - f_1(x_1) + v_1, \quad x_i^+ = x_i - f_i(x_i) + v_i + (1 - P_{i-1})f_{i-1}(x_{i-1}), \text{ for } i = 2,...,n. \quad (3.22)$$

In what follows, we assume that $x^* = (x_1^*,...,x_n^*) \in \prod_{i=1}^{n}(0, \mu_i - \varepsilon]$, $v_i^* \in [b_i + \varepsilon, v_{i,\max}]$ for $i \in R$ and for some $\varepsilon \in (0, 1/2)$ and $v^* \in (0, v_{1,\max}] \times \prod_{i=2}^{n}[0, v_{i,\max}]$. Moreover, we assume that $P_i \in [0, 1-\varepsilon]$ for $i = 1,...,n-1$ and $r_i \in [\varepsilon, 1-\varepsilon]$ for $i = 1,...,n$.

Another feature of the present problem is that the selection of the uncongested equilibrium point may be made in an implicit way. For example, we may want the uncongested equilibrium point that guarantees the maximum outflow from the freeway. In such cases, the equilibrium position of the controllable inflows is determined as a function of the nominal values of the uncontrollable inflows and the parameters of the freeway, i.e., there exists a smooth function

$$g: [0, 1-\varepsilon]^{n-1} \times \prod_{i \notin R}[0, v_{i,\max}] \times [\varepsilon, 1-\varepsilon]^n \to \prod_{i \in R}[b_i + \varepsilon, v_{i,\max}]$$

such that

$$(v_i^*; i \in R) = g(P, v_i^*; i \notin R, r) \quad (3.23)$$

where $P = (P_1,...,P_{n-1})' \in [0, 1-\varepsilon]^{n-1}$ and $r = (r_1,...,r_n)' \in [\varepsilon, 1-\varepsilon]^n$.

## 3.III. Measurements and Unknown Parameters

Let $m \in \{1,...,n\}$ be the cardinal number of the set $R$ and let $u \in U = \prod_{i \in R}[b_i, v_{i,\max}] \subseteq (0,+\infty)^m$ be the vector of all controllable inflows $v_i$ with $i \in R$.

The model parameters which are (usually) unknown or uncertain are: the exit rates $P_i \in [0,1)$ for $i = 1,...,n-1$, the uncontrollable inflows $v_i^* \in \Re_+$ for $i \notin R$ and the demand coefficients $r_i \in (0,1)$ for $i = 1,...,n$. All these parameters will be denoted by $\theta = (P, v_i^*; i \notin R, r)$ and are assumed to take values in a compact set $\Theta := [0, 1-\varepsilon]^{n-1} \times \prod_{i \notin R}[0, v_{i,\max}] \times [\varepsilon, 1-\varepsilon]^n$, for some $\varepsilon \in (0, 1/2)$. Therefore, the control system (3.11), (3.12), (3.13) can be written in the following vector form:

$$\begin{aligned} x^+ &= \overline{F}(d, \theta, x, u) \\ x \in S, d \in D, \theta &\in \Theta, u \in U = \prod_{i \in R}[b_i, v_{i,\max}] \end{aligned} \quad (3.24)$$

Notice that the feedback law defined by (3.18) is a feedback law of the form $u = k(\theta, x)$: the feedback law depends on the unknown parameters through $x^*$ and $(v_i^*; i \in R)$ (recall (3.17) and (3.23)). It follows that assumption (H1) holds for system (3.24). An explicit definition of the feedback law $k: \Theta \times S \to U$ is given by the following equations for all $\hat{\theta} = (\hat{P}, \hat{v}_i^*; i \notin R, \hat{r}) \in \Theta$, $x \in S$ with $\hat{r} = (\hat{r}_1,...,\hat{r}_n)' \in [\varepsilon, 1-\varepsilon]^n$, $\hat{P} = (\hat{P}_1,..., \hat{P}_{n-1})' \in [0, 1-\varepsilon]^{n-1}$:

$$(\hat{v}_i^*; i \in R) = g(\hat{P}, \hat{v}_i^*; i \notin R, \hat{r}) \quad (3.25)$$

$$\begin{aligned} \hat{x}_1^* &= \min\left(\hat{r}_1^{-1} \hat{v}_1^*, \mu_1 - \varepsilon\right) \\ \hat{x}_i^* &= \min\left(\hat{r}_i^{-1}\left(\hat{v}_i^* + \sum_{j=1}^{i-1} \hat{v}_j^* \left(\prod_{k=j}^{i-1}(1-\hat{P}_k)\right)\right), \mu_i - \varepsilon\right) \quad , \quad i = 2,...,n \end{aligned} \quad (3.26)$$

$$u = k(\hat{\theta}, x) \text{ with } k_i(\hat{\theta}, x) = \max\left(b_i, \hat{v}_i^* - \tau^{-1}\left(\hat{v}_i^* - b_i\right)\Xi(\hat{\theta}, x)\right), \text{ for all } x \in S, i \in R \quad (3.27)$$

$$\Xi(\hat{\theta}, x) := \sum_{i=1}^n \sigma^i \max\left(0, x_i - \hat{x}_i^*\right), \text{ for all } x \in S, \quad (3.28)$$

The measured quantities are the cell densities $x \in S$ and the outflows from each cell. We have two kinds of outflows from each cell: the outflow to regions out of the freeway

$$\begin{aligned} Q_{out} &= (Q_{1,out},...,Q_{n,out})' \in \Re_+^n \\ Q_{i,out} &= P_i s_{i+1} f_i(x_i) \quad , \quad i = 1,...,n-1 \\ Q_{n,out} &= f_n(x_n) \end{aligned} \quad (3.29)$$

and the outflows from one cell to the next cell

$$\begin{aligned} Q &= (Q_1,...,Q_{n-1})' \in \Re_+^{n-1} \\ Q_i &= (1-P_i) s_{i+1} f_i(x_i) \quad , \quad i = 1,...,n-1 \end{aligned} \quad (3.30)$$

Therefore, the measured output is given by:

$$y = h(d, \theta, x) = (x, Q_{out}, Q) \in S \times \Re_+^n \times \Re_+^{n-1}, \quad (3.31)$$

Assumption (H) guarantees that $h(D \times \Theta \times S) \subseteq Y := S \times \prod_{i=1}^n [0, a_i] \times \prod_{i=1}^{n-1}[0, a_i]$. Notice that $Y := S \times \prod_{i=1}^n [0, a_i] \times \prod_{i=1}^{n-1}[0, a_i]$ is a bounded set.

It follows from (3.21), (3.22), (3.29), (3.30), assumption (H) and the fact that $\mu_i \in (0, \delta_i)$ ($i = 1,...,n$), that:

"if $x(t-1) \in \Omega = \prod_{i=1}^{n}(0,\mu_i)$, $t \geq 1$, then the following equations hold:

$$P_i = \frac{Q_{i,out}(t-1)}{Q_{i,out}(t-1)+Q_i(t-1)}, \quad i=1,...,n-1 \tag{3.32}$$

$$v_i^* = x_i(t) - x_i(t-1) + Q_i(t-1) + Q_{i,out}(t-1) - Q_{i-1}(t-1), \quad i \in \{2,...,n\}\setminus R \tag{3.33}$$

$$v_1^* = x_1(t) - x_1(t-1) + Q_1(t-1) + Q_{1,out}(t-1), \text{ if } 1 \notin R \tag{3.34}$$

$$r_i = \frac{Q_{i,out}(t-1)+Q_i(t-1)}{x_i(t-1)}, \quad i=1,...,n". \tag{3.35}$$

Equations (3.32), (3.33), (3.34), (3.35), (3.31) allow us to define a mapping $\Psi: h(D\times\Theta\times S)\times Y \to \Theta$ for which $\theta = (P_1,...,P_{n-1}, v_i^*; i \notin R, r_1,...,r_n)' = \Psi(y(t), y(t-1))$ for all $t \geq 1$ with $y(t-1) \in A$, where $A \subseteq Y$ is the set for which

$$w=(w_1,w_2,w_3) \in A \Leftrightarrow (w_1,w_2,w_3) \in Y, \ w_1 \in \Omega = \prod_{i=1}^n(0,\mu_i) \text{ and } w_{2,i}+w_{3,i} > 0 \text{ for } i=1,...,n-1. \tag{3.36}$$

The mapping $\Psi: h(D\times\Theta\times S)\times Y \to \Theta$ is defined by

$$\hat{\theta} = (\hat{P}_1,...,\hat{P}_{n-1}, \hat{v}_i^*; i \notin R, \hat{r}_1,...,\hat{r}_n)' = \Psi(y,w) \tag{3.37}$$

$$\hat{P}_i = \min\left(1-\varepsilon, \frac{w_{2,i}}{w_{2,i}+w_{3,i}}\right), \quad i=1,...,n-1 \tag{3.38}$$

$$\hat{v}_i^* = \max(0, \min(v_{i,\max}, x_i - w_{1,i} + w_{3,i} + w_{2,i} - w_{3,i-1})), \text{ if } i \in \{2,...,n\}\setminus R \text{ and } i \neq n \tag{3.39}$$

$$\hat{v}_n^* = \max(0, \min(v_{n,\max}, x_n - w_{1,n} + w_{2,n} - w_{3,n-1})), \text{ if } n \notin R \tag{3.40}$$

$$\hat{v}_1^* = \max(0, \min(v_{1,\max}, x_1 - w_{1,1} + w_{3,1} + w_{2,1})), \text{ if } 1 \notin R \tag{3.41}$$

$$\hat{r}_i = \max\left(\varepsilon, \min\left(1-\varepsilon, \frac{w_{2,i}+w_{3,i}}{w_{1,i}}\right)\right), \quad i=1,...,n-1 \tag{3.42}$$

$$\hat{r}_n = \max\left(\varepsilon, \min\left(1-\varepsilon, \frac{w_{2,n}}{w_{1,n}}\right)\right) \tag{3.43}$$

Using assumption (H), (3.16), (3.17) and (3.31), it follows that there exists $y^* \in Y$ with $y^* = h(d, \theta, x^*)$ for all $d \in D$. By virtue of our assumption $x^* = (x_1^*,...,x_n^*) \in \prod_{i=1}^n(0,\mu_i)$ and $v^* \in (0, v_{1,\max}]\times\prod_{i=2}^n[0, v_{i,\max}]$, (3.36), we conclude that $A$ contains all $w \in Y$ in a neighborhood of $y^*$. It follows that (H2) holds with $p=1$ for system (3.24) with output given by (3.29), (3.30), (3.31).

In order to prove that assumption (H3) holds for system (3.24) with output given by (3.29), (3.30), (3.31), we need the following fact, which is a consequence of property (C5) shown in [14] and (3.20).

**Fact:** *Define* $I_j(x) := \sum_{i=1}^{j} x_i$ *for* $j=1,...,n$. *There exists a constant* $C \in (0,1)$ *such that the following inequality holds:*

$$\sum_{i=1}^n I_i(x^+) \leq (1-C)\sum_{i=1}^n I_i(x) + \sum_{i=1}^n (n+1-i)v_i, \text{ for all } (x,v,d) \in S\times(0,v_{1,\max}]\times\prod_{i=2}^n[0,v_{i,\max}]\times[0,1]^{n-1} \tag{3.44}$$

*where $x^+$ is given by (3.24).*

The following proposition guarantees that assumption (H3) holds for system (3.24) with output (3.29), (3.30), (3.31).

**Proposition 3.2:** *Suppose that $b_i > 0$ ($i \in R$) and $v_{i,\max}$ ($i \notin R$) are sufficiently small and that $\tau > 0$ is sufficiently small ($\tau \leq \varepsilon^2 \sigma^n \min_{i \in R}\left((v_{i,\max} - b_i)^{-1}\right)$). Then there exists an integer $m \geq 1$ such that for every sequence $\{(d(t), \hat{\theta}(t)) \in D \times \Theta\}_{t=0}^{\infty}$ and for every $x_0 \in S$, the solution $x(t)$ of (3.24), (3.31) with $u = k(\hat{\theta}, x)$, initial condition $x(0) = x_0$ corresponding to inputs $\{(d(t), \hat{\theta}(t)) \in D \times \Theta\}_{t=0}^{\infty}$ satisfies $y(t - 1 - i(t)) \in A$ for some $i(t) \in \{0, 1, ..., m\}$ and for all $t \geq m + 1$.*

**Proof:** Assume that $b_i > 0$ ($i \in R$) and $v_{i,\max}$ ($i \notin R$) are sufficiently small so that

$$\sum_{i \in R}(n+1-i)b_i + \sum_{i \notin R}(n+1-i)v_{i,\max} < C \min_{i=1,...,n}((n+1-i)\mu_i). \tag{3.45}$$

Since $\tau \leq \varepsilon^2 \sigma^n \min_{i \in R}\left((v_{i,\max} - b_i)^{-1}\right)$ and $\hat{v}_i^* \in [b_i + \varepsilon, v_{i,\max}]$ for $i \in R$, it follows that

$$\tau^{-1}(\hat{v}_i^* - b_i) \geq \varepsilon^{-1}(v_{i,\max} - b_i)\sigma^{-n}, \text{ for all } i \in R. \tag{3.46}$$

Let $m \geq 1$ be an integer that satisfies

$$m \geq 2 + \left\lceil \frac{\ln\left(\min_{i=1,...,n}((n+1-i)\mu_i) - C^{-1}\kappa\right) - \ln\left(\sum_{i=1}^{n}(n+1-i)a_i\right)}{\ln(1-C)} \right\rceil. \tag{3.47}$$

Next, we show the following claim.

**Claim:** If $x \notin \Omega$ then for every $(\hat{\theta}, d) \in \Theta \times [0,1]^{n-1}$ it holds that

$$\sum_{i=1}^{n} I_i(x^+) \leq (1-C)\sum_{i=1}^{n} I_i(x) + \kappa \tag{3.48}$$

where $C \in (0,1)$ is the constant involved in (3.44), $\kappa := \sum_{i \in R}(n+1-i)b_i + \sum_{i \notin R}(n+1-i)v_{i,\max}$ and $x^+$ is given by (3.24) with $u = k(\hat{\theta}, x)$.

**Proof of Claim:** If $x \notin \Omega = \prod_{i=1}^{n}(0, \mu_i)$, then there exists $i \in \{1,...,n\}$ such that $x_i \geq \mu_i$. Since $\hat{x}^* = (\hat{x}_1^*, ..., \hat{x}_n^*) \in \prod_{i=1}^{n}[0, \mu_i - \varepsilon]$ (recall (3.26)), it follows from (3.28) and the fact that $\sigma \in (0,1]$ that $\Xi(\hat{\theta}, x) \geq \sigma^n(x_i - \hat{x}_i^*) \geq \varepsilon \sigma^n$. Since (3.46) holds, it follows from (3.27) that $v_i = b_i$ for all $i \in R$. Inequality (3.48) is a consequence of (3.44) and the fact that $v_i^* \in [0, v_{i,\max}]$ for all $i \notin R$. The proof of the claim is complete.

We show next, by means of a contradiction, that for every sequence $\{(d(t), \hat{\theta}(t)) \in D \times \Theta\}_{t=0}^{\infty}$ and for every $x_0 \in S$, the solution $x(t)$ of (3.24), (3.31) with $u = k(\hat{\theta}, y)$, initial condition $x(0) = x_0$ corresponding to inputs $\{(d(t), \hat{\theta}(t)) \in D \times \Theta\}_{t=0}^{\infty}$ satisfies $y(t - 1 - i(t)) \in A$ for some $i(t) \in \{0,1,...,m\}$ and for all $t \geq m + 1$.

Suppose that, on the contrary, there exists a sequence $\{(d(t), \hat{\theta}(t)) \in D \times \Theta\}_{t=0}^{\infty}$, a vector $x_0 \in S$ and an integer $t \geq m + 1$, such that the solution $x(t)$ of (3.24), (3.31) with $u = k(\hat{\theta}, y)$, initial condition $x(0) = x_0$ corresponding to inputs $\{(d(t), \hat{\theta}(t)) \in D \times \Theta\}_{t=0}^{\infty}$ satisfies $y(t - 1 - i(t)) \notin A$ for all $i(t) \in \{0,1,...,m\}$. By virtue of (3.36), this implies that $x(t - 1 - i(t)) \notin \Omega$ for all $i(t) \in \{0,1,...,m\}$ (notice that (3.21), (3.22), (3.29), (3.30), (3.31) and (3.36) guarantee that $x \in \Omega$ implies that $y \in A$). It follows from the Claim, that

$$\sum_{i=1}^{n} I_i(x(l+1)) \le (1-C) \sum_{i=1}^{n} I_i(x(l)) + \kappa, \text{ for } l = t-1-m,\ldots,t-1. \tag{3.49}$$

Using (3.49) repeatedly, we get:

$$\sum_{i=1}^{n} I_i(x(t-1)) \le (1-C)^m \sum_{i=1}^{n} I_i(x(t-1-m)) + \kappa \frac{1-(1-C)^m}{C}. \tag{3.50}$$

Using the definition $I_j(x) := \sum_{i=1}^{j} x_i$ for $j = 1,\ldots,n$ and the fact that $x \in S = \prod_{i=1}^{n}(0, a_i]$, we get from (3.50):

$$(n+1-j)x_j(t-1) \le (1-C)^m \sum_{i=1}^{n}(n+1-i)a_i + C^{-1}\kappa, \text{ for all } j = 1,\ldots,n. \tag{3.51}$$

Using (3.51), (3.45) and (3.47), we get:

$$(n+1-j)x_j(t-1) < \min_{i=1,\ldots,n}((n+1-i)\mu_i), \text{ for all } j = 1,\ldots,n$$

which implies that $x(t-1) \in \Omega = \prod_{i=1}^{n}(0,\mu_i)$, a contradiction. The proof is complete. ◁

The main result for the freeway model is a consequence of Theorem 2.1 and the fact that all functions are sufficiently smooth in a neighborhood of the equilibrium.

**Corollary 3.3:** *Consider system (3.24) with output given by (3.29), (3.30), (3.31). Suppose that $b_i > 0$ $(i \in R)$ and $v_{i,\max}$ $(i \notin R)$ are sufficiently small and that $\tau > 0$ is sufficiently small. Then the dynamic feedback law given by:*

$$w_1^+ = x, \quad w_2^+ = Q_{out}, \quad w_3^+ = Q \tag{3.52}$$

$$\hat{P}_i^+ = \begin{cases} \hat{P}_i & \text{if } w \notin A \\ \min\left(1-\varepsilon, \dfrac{w_{2,i}}{w_{2,i}+w_{3,i}}\right) & \text{if } w \in A \end{cases}, \quad i = 1,\ldots,n-1 \tag{3.53}$$

$$(\hat{v}_i^*)^+ = \begin{cases} \hat{v}_i^* & \text{if } w \notin A \\ \max(0, \min(v_{i,\max}, x_i - w_{1,i} + w_{3,i} + w_{2,i} - w_{3,i-1})) & \text{if } w \in A \end{cases}, \text{ if } i \in \{2,\ldots,n\} \setminus R \text{ and } i \ne n \tag{3.54}$$

$$(\hat{v}_n^*)^+ = \begin{cases} \hat{v}_n^* & \text{if } w \notin A \\ \max(0, \min(v_{n,\max}, x_n - w_{1,n} + w_{2,n} - w_{3,n-1})) & \text{if } w \in A \end{cases}, \text{ if } n \notin R \tag{3.55}$$

$$(\hat{v}_1^*)^+ = \begin{cases} \hat{v}_1^* & \text{if } w \notin A \\ \max(0, \min(v_{1,\max}, x_1 - w_{1,1} + w_{3,1} + w_{2,1})) & \text{if } w \in A \end{cases}, \text{ if } 1 \notin R \tag{3.56}$$

$$\hat{r}_i^+ = \begin{cases} \hat{r}_i & \text{if } w \notin A \\ \max\left(\varepsilon, \min\left(1-\varepsilon, \dfrac{w_{2,i}+w_{3,i}}{w_{1,i}}\right)\right) & \text{if } w \in A \end{cases}, \quad i = 1,\ldots,n-1 \tag{3.57}$$

$$\hat{r}_n^+ = \begin{cases} \hat{r}_n & \text{if } w \notin A \\ \max\left(\varepsilon, \min\left(1-\varepsilon, \dfrac{w_{2,n}}{w_{1,n}}\right)\right) & \text{if } w \in A \end{cases} \tag{3.58}$$

*with (3.25), (3.26), (3.27), (3.28), $\hat{P} = (\hat{P}_1,\ldots,\hat{P}_{n-1})$, $P = (P_1,\ldots,P_{n-1})$, $\hat{r} = (\hat{r}_1,\ldots,\hat{r}_n)$, $r = (r_1,\ldots,r_n)$, $w = (w_1, w_2, w_3)$, $\hat{v}^* = (\hat{v}_1^*,\ldots,\hat{v}_n^*)$, achieves the following:*

*1) There exist constants $M, \sigma > 0$ such that for every sequence $\{d(i) \in D\}_{i=0}^{\infty}$ and for every $(x_0, w_0, \hat{P}_0, \hat{v}_j^*; j \notin R, \hat{r}_0) \in S \times Y \times \Theta$, the solution of the closed-loop system (3.24), (3.31) with (3.52)-(3.58), (3.25)-(3.28), initial condition $(x(0), w(0), \hat{p}(0), \hat{v}_j^*(0); j \notin R, \hat{r}(0)) = (x_0, w_0, \hat{p}_0, \hat{v}_j^*; j \notin R, \hat{r}_0)$ corresponding to input $\{d(i) \in D\}_{i=0}^{\infty}$ satisfies*

$$|x(t)-x^*|+|w(t)-y^*|+|\hat{r}(t)-r|+|\hat{P}(t)-P|+|\hat{v}^*(t)-v^*|$$

$$\leq M\exp(-\sigma t)\left(|x(0)-x^*|+|w(0)-y^*|+|\hat{r}(0)-r|+|\hat{P}(0)-P|+\sum_{i\notin R}|\hat{v}_i^*-v_i^*|\right), \text{ for all } t\geq 0 \quad (3.59)$$

2) There exists an integer $N\geq 1$ such that for every sequence $\{d(i)\in D\}_{i=0}^{\infty}$ and for every $(x_0,w_0,\hat{P}_0,\hat{v}_j^*;j\notin R,\hat{r}_0)\in S\times Y\times\Theta$, the solution of the closed-loop system (3.24), (3.31) with (3.52)-(3.58), (3.25)-(3.28), initial condition $(x(0),w(0),\hat{P}(0),\hat{v}_j^*(0);j\notin R,\hat{r}(0))=(x_0,w_0,\hat{P}_0,\hat{v}_j^*;j\notin R,\hat{r}_0)$ corresponding to input $\{d(i)\in D\}_{i=0}^{\infty}$ satisfies $\hat{P}(t)=P$, $\hat{r}(t)=r$, $\hat{v}^*(t)=v^*$, for all $t\geq N$.

**Proof:** Let $N_1\subseteq\Omega$ be a neighborhood of $x^*$, $N_2\subseteq A$ be a neighborhood of $y^*$, and let $N_3\subseteq\Re^{3n-1-m}$ be a neighborhood of $\theta$. Since $\Omega=\prod_{i=1}^{n}(0,\mu_i)$, it follows from Assumption (H) and the fact that $\mu_i\in(0,\delta_i)$ for $i=1,...,n$ that $f_i(x_i)=r_ix_i$ for $i=1,...,n$. Definitions (3.29), (3.30), (3.31) in conjunction with (3.21) and the fact that $P_i\in[0,1)$ for $i=1,...,n-1$, $r_i\in(0,1)$ for $i=1,...,n$, imply that the following inequality holds for all $x\in\Omega=\prod_{i=1}^{n}(0,\mu_i)$ and $d=(d_2,...,d_n)\in D=[0,1]^{n-1}$:

$$\begin{aligned}|h(d,\theta,x)-y^*| &\leq |x-x^*|+|Q_{out}-Q_{out}^*|+|Q-Q^*| \\ &\leq |x-x^*|+\sum_{i=1}^{n-1}|P_if_i(x_i)-P_if_i(x_i^*)|+|f_n(x_n)-f_n(x_n^*)| \\ &\quad +\sum_{i=1}^{n-1}|(1-P_i)f_i(x_i)-(1-P_i)f_i(x_i^*)| \\ &\leq |x-x^*|+\sum_{i=1}^{n}r_i|x_i-x_i^*|\leq\left(1+\sum_{i=1}^{n}r_i\right)|x-x^*|\end{aligned} \quad (3.60)$$

Next, we notice that by virtue of (3.22) and the facts that $P_i\in[0,1)$ for $i=1,...,n-1$, $r_i\in(0,1)$ for $i=1,...,n$, $f_i(x_i^*)=v_i^*+(1-P_{i-1})f_{i-1}(x_{i-1}^*)$ for $i=2,...,n$, $f_1(x_1^*)=v_1^*$, it follows that the following holds for all $x\in\Omega=\prod_{i=1}^{n}(0,\mu_i)$, $d\in D=[0,1]^{n-1}$ and $u\in\Re^m$:

$$\begin{aligned}&|\overline{F}(d,\theta,x,u)-x^*| \\ &\leq |x_1-f_1(x_1)+v_1-x_1^*|+\sum_{i=2}^{n}|x_i-f_i(x_i)+v_i+(1-P_{i-1})f_{i-1}(x_{i-1})-x_i^*| \\ &\leq \sum_{i=2}^{n}|x_i-f_i(x_i)+f_i(x_i^*)+(1-P_{i-1})f_{i-1}(x_{i-1})-(1-P_{i-1})f_{i-1}(x_{i-1}^*)-x_i^*| \\ &\quad +m|u-u^*|+|x_1-f_1(x_1)+f_1(x_1^*)-x_1^*| \\ &\leq (1-r_1)|x_1-x_1^*|+\sum_{i=2}^{n}(1-r_i)|x_i-x_i^*| \\ &\quad +\sum_{i=2}^{n}(1-P_{i-1})r_{i-1}|x_{i-1}-x_{i-1}^*|+m|u-u^*| \\ &\leq \left(n-\sum_{i=1}^{n}r_i+\sum_{i=2}^{n}(1-P_{i-1})r_{i-1}\right)|x-x^*|+m|u-u^*|\end{aligned} \quad (3.61)$$

where $u^*=(v_i^*;i\in R)$. Using (3.27) and (3.28), it is straightforward to show that there exists a constant $\tilde{L}>0$ such that the following inequality holds for all $x,\hat{x}^*\in S$ and $\hat{v}^*\in\prod_{i=1}^{n}[0,v_{i,\max}]$:

$$|u-u^*|\leq\tilde{L}|x-x^*|+\tilde{L}|\hat{x}^*-x^*|+\tilde{L}|\hat{v}^*-v^*|. \quad (3.62)$$

Using (3.25), (3.26) and the fact that the function $g:[0,1-\varepsilon]^{n-1} \times \prod_{i \notin R}[0,v_{i,\max}] \times [\varepsilon,1-\varepsilon]^n \to \prod_{i \in R}[b_i+\varepsilon,v_{i,\max}]$ is a smooth function, it follows that the following inequality holds for all $\hat{\theta} \in N_3 \cap \Theta$:

$$\left|\hat{x}^* - x^*\right| + \left|\hat{v}^* - v^*\right| \leq M\left|\hat{\theta} - \theta\right|. \tag{3.63}$$

Finally, using definitions (3.37)-(3.43) in conjunction with the fact that $N_2 \subseteq A$, it follows that there exists a constant $\bar{L} > 0$ such that

$$\left|\Psi(h(d,\theta,x),w) - \theta\right| \leq \bar{L}\left|x - x^*\right| + \bar{L}\sum_{i=1}^{p}\left|w_i - y^*\right|,$$

for all $x \in N_1 \cap S$, $d \in D$, $\hat{\theta} \in N_3 \cap \Theta$, $w_i \in N_2 \cap Y$ ($i=1,...,p$) with $w = (w_1,...,w_p)$. (3.64)

Since, we have already proved that assumptions (H1), (H2), (H3) hold for the closed-loop system (3.24), (3.31) with (3.52)-(3.58), (3.25)-(3.28), it follows from (3.60), (3.61), (3.62), (3.63) and (3.64) that all assumptions of Theorem 2.1 hold. Therefore, Corollary 3.3 is a direct application of Theorem 2.1 to the closed-loop system (3.24), (3.31) with (3.52)-(3.58), (3.25)-(3.28). The proof is complete. ◁

## 4. Concluding Remarks

Novel results for adaptive control schemes for uncertain discrete-time systems, which guarantee robust, global, exponential convergence to the desired equilibrium point of the system, were provided in the present work. The proposed control scheme consists of a nominal feedback law, which achieves robust, global, exponential stability properties when the vector of the parameters is known, in conjunction with a nonlinear, dead-beat observer. The proposed adaptive scheme did not require the knowledge of a Lyapunov function for the closed-loop system under the action of the nominal feedback stabilizer and is directly applicable to highly nonlinear, uncertain discrete-time systems with unknown constant parameters.

**Acknowledgments:** The research leading to these results has received funding from the European Research Council under the E.U.'s 7[th] Framework Programme (FP/2007-2013)/ERC Grant Agreement n. [321132], project TRAMAN21.